\documentclass[11pt]{amsart}

\usepackage{amssymb, amsmath}

\usepackage[polutonikogreek,english]{babel}
\usepackage[GlyphNames]{teubner}
\usepackage{geometry}
\usepackage{verbatim}

\begin{document}
\title[Mathematical Interpretation  of Plato's Third Man Argument]
{Mathematical Interpretation  of Plato's Third Man Argument based on the notion of Convergence}
\author {George Chailos}
 \address{Department of Computer Science, University of Nicosia , Nicosia 1700, Cyprus }
\email{chailos.g@unic.ac.cy}

\subjclass[2010]{62A01, 97E20, 03A05, 97E30}
\keywords{Platonic Philosophy, Convergence, Complete Continuum, Infinity.}

\begin {abstract}
 The main aim of this article is to defend the thesis that Plato apprehended the structure of incommensurable magnitudes in a way that these magnitudes correspond in a unique and well defined manner to the modern concept of the \textsl{Dedekind cut}. Thus, the notion of convergence is
 consistent with Plato's apprehension of mathematical concepts, and in particular these of  \textsl{density} of magnitudes and
 the \textsl{complete continuum} in the sense that they include incommensurable cuts.
 For this purpose I discuss and interpret, in a new perspective,  the mathematical framework
 and the logic of the Third Man Argument (TMA) that appears in Plato's  \textit{Parmenides} as well as mathematical concepts from other Platonic dialogues.
 I claim that in this perspective the apparent infinite sequence of $F$-Forms, that it is generated
 by repetitive applications of the TMA, converges (in a mathematical sense) to a unique $F$-Form for the particular predicate.
 I also claim and prove that within this framework the logic of the TMA is consistent
 with that of the Third Bed Argument (TBA) as presented in Plato's \textit{Republic}.
 This supports Plato's intention for assuming a unique Form per Predicate; that is, the Uniqueness thesis.
\end{abstract}

\maketitle

\section{Introduction}

 In this article I aim to provide an adequate mathematical interpretation of the classical Third Man
Argument (TMA) which is strongly related with the mathematical approach of Plato's theory of Forms.

In this section I introduce the problem, and in the second section I present
the framework within  which my arguments are based on, and the Main Claim is founded. In the third section I
analyze, defend and eventually provide an adequate proof of the Main Claim. In the fourth section  I present in detail the
mathematical concepts and  the topological framework within which my arguments and the main thesis are comprehended and analyzed.
In the last section, summarizing the work done in the previous sections, I briefly present the main conclusions of my work.

The literature dealing with TMA, already large in $1954$, has become enormous since then and almost all of the authors have followed
Gregory Vlastos, where in his famous paper of 1954 \cite{vlastos1} pointed out  that the argument is formally a \textit{non sequitur}.

The classical TMA appears in the Plato's dialogue
\textit{Parmenides} $132a1-b2$\footnote{Text: `\textsl{...This, I suppose, is what leads you to believe that each form is one.
Whenever many things seem to you to be large, some one form probably seems to you to be the same when you look at them all. So you think
that largeness is one. . . . But what about largeness itself and the other large things? If you look at them all in your mind in the same
way, won't some one largeness appear once again, by virtue of which they all appear large? . . . So another form of largeness will have
made an appearance, besides largeness itself and its participants. And there will be yet another over all these, by virtue of which
they will all be large. So each of your forms will no longer be one, but an infinite multitude}...'.
(The translation is taken from Cohen and Keyt \cite{cohen3}).}, and essentially its logic
is present in other arguments in Platonic dialogues, such as the TMA version in
\textit{Parm.}\,$132d1-133a6$, where it is applied on the `Form-idea' of `Resemblance'. I aim to adequately explain the
apparent plurality of Forms (for a certain predicate) that appears in the TMA and to show that it is compatible with Plato's thesis
for the existence of a single form per predicate; that is the uniqueness thesis\footnote{Plato's intention in defending the
uniqueness of a Form per predicate was clearly introduced in the \textit{Third Bed Argument}, TBA, \textit{Republic} $597c-d$,
and is also present in other Platonic texts. The related phrase `\textgreek{\era n{} \era kaston{} e\isc dos}' in \textit{Parm.} $132a1$
and its relation to the uniqueness thesis is analyzed by Cohen \cite{cohen1}, pp.$433-466$.
The uniqueness thesis  should be expressed as:\\
`\textsl{There is exactly one Form corresponding to every predicate that has a Form}'.\\
For a different version of the uniqueness thesis and an extensive discussion the reader should
consult G. Fine \cite{fine}. \label{fnt:unique}}.
For supporting this, I take into account the Platonic apprehension of fundamental mathematical concepts,
the Mathematics developed in Plato's Academy as well as  Plato's dialectic.

The Theory of Forms is also a theory of judgment.
Judging involves consulting Forms: To judge that an object $x$, \emph{either a sensible particular or a Form}, is $F$,
is to consult the form of $F$-ness and to perceive $x$ as being sufficiently
like $F$-ness to qualify for the predicate $F$. Alan Code in \cite{code} suggests that the TMA raises an objection to this theory of judgment.
I quote:

\begin{quote}
 {\small`\textsl{The TMA is designed to reduce to absurdity the claim that it is the consultation of Forms which enables us to
 make judgments. It does this by showing that if that were the case, we would have to perform an infinite number of such
 consultations to make just one judgment}.'}
\end{quote}

Since the theory of Forms tries to \emph{explain} predication, the TMA is also a challenge to it as a theory of predication.
It is evident though that `participating in a Form' is supposed to explain predication. And the upshot of the TMA, as was presented
by many authors, is that there is something defective about this explanation: Since trying to explain predication
in terms of the notion of participating in a paradigmatic Form leads to an infinite regress, and hence is no explanation at all.


In developing and defending my arguments, I shall  be consistent with the interpretation of the \emph{presence of a property in a thing},
as well as the \emph{recurrence of a single property} in different things.
According to Scaltsas \cite{scaltsas}, the things are $F$ by participating in a Form $F$-ness is the answer to two different
questions that Plato implies in \textit{Phaedo} $100c9-d8$, \textit{Parm.} $128e6-129a4$, $130e5-131a2$, as well as earlier
in \textit{Meno} $72c$. The first is `Why is a thing $F$?' That is, the first question concerns the predication of $F$-ness. In
\textit{Phaedo} $100c9-d8$ the forms are introduced as the causes of things being $F$\footnote{In \textit{Phaedo} \emph{the} Form
is also referred as \textsl{the $F$ itself} ($74a11-12$), which is the
\textgreek{kaj> a\r ut\og} $F$ ($74b3-4$), \textsl{the cause} (\textgreek{a\is t\ia a}) that \textsl{makes}
(\textgreek{poie\ic}) things being $F$ ($100c9-d8)$, or the explanation of something being $F$. \label{fnt:itself}}.
The second question is `Why are different things similar?' This  question that appears clearly
in \textit{Parm.} $128e6-129a4, 130e5-131a2$ concerns recurrence of a single property in different things
and considers the quality identity with respect to $F$-ness.

In this article I defend the thesis  that the apparent infinite sequence
of Forms, $\{F_i\}_{i=0}^\infty$, $i \in \mathbb{N}$\footnote{where $\mathbb{N}$ is the set of Natural Numbers.},
that appears by  repetitive applications of the TMA,
for either explaining predication or making a judgment that $x$ is $F$, it is
\emph{increasing} (ascending) and has a (mathematical) \emph{limit};
that is, an attainable \textsl{least upper bound}. In other words, this infinite sequence $\{F_i\}_{i=0}^\infty$
converges to a unique Form $F$ for a particular predicate.
Henceforth, $\displaystyle{\underset{i \to +\infty} \lim F_i=F}$,  where the convergence is understood as a mathematical
one\footnote{The formal definition of the convergence using the concept of a limit is: $\forall \epsilon>0$, there is a natural number $k_{0}$,
such that $\||{F-F_i}\||<\epsilon$, whenever $n\geq k_{0}$. The $\||\cdot\||$ denotes the norm of the vector-topological space in concern.
See also Rudin \cite{rudin} for the mathematical definitions of the \textit{limit} and the \textit{least upper bound} (\textit{l.u.b}.).}.
This terminating-limit $F$- Form should be also apprehended as compatible to the
`\textsf{anupotheton arx\={e}n}' -`\textgreek{\as nup\oa jeton{} \as rq\hg n}'- in \textit{Republic's}
language\footnote{The term appears in \textit{Rep.} $509b-511d$, in $510b$ and in $511b$,
and shall be interpreted later in the paper. We do not give any translation, since any translation may lead
to a specific interpretation. (In \cite{thesaurus} \textit{ad.loc.} is translated as `the principle that
transcends assumption'.) For an extensive analysis of this passage and especially the concept and the status
of this term we refer to  Karasmanis \cite{karasm2}, \cite{karasm3} and  Benson \cite{benson}.}, but applied here for
each particular predicate $F$. Moreover, this $F$-Form  should be also considered as analogous to the
`\textsl{final rung of Diotima's ladder}' as presented in \textit{Symposium} $210e$\textit{ff}\footnote{We shall see that
this approach is compatible with the mathematics of Plato's academy and Plato's dialectic.
For an elaborate and comprehensive exposition of Plato's dialectic see Robinson \cite{robinson} ch. $6,7$ and $10$.}.
(In \textit{Symp.} the above procedure is developed in the context of a particular Form, namely the Form of `Beauty'.)

Furthermore, I shall argue and justify that the logics of the TMA  and the TBA arguments do not contradict each other;
instead they are consistent and mutually complementary.
In this perspective it is necessary to apprehend this terminating $F$-Form, that is obtained as a (mathematical) limit of
the above infinite sequence, in a somewhat different context than that of the various $\{F_{i}\}_{i=0}^{\infty}$ , $i\in \mathbb{N}$,
in the sequence. This thesis is developed, properly analyzed and defended in the sequel.

 \section{Foundation of the Argument}
 Before proceeding, I shall briefly present the necessary historical background of the problem for developing
and defending my arguments and my thesis.

In the following, the schematic letter `$F$' shall serve as a dummy predicate for any predicate for which there is a Form.
(It is typically used in place of predicate `large' that appears in the TMA).
We also note that we shall not deal with any issues related to the so called `Imperfection Argument'\footnote{ For example this argument
does not posit a form even for every property-name;
it posits a form for the predicate \textit{large} but not a Form for the predicate \textit{man}.
And it supports that we can infer that is a Form of $F$ only when we have a group that consists of imperfectly $F$ things.
Namely, the imperfection argument posits Forms  both for restricted range of predicates and also a restricted range of groups.
(For further details on the `Imperfection Argument' we refer to G. Fine's \cite{fine}.).}
(as entailed primarily in \textit{Rep.} $523-525$, or elsewhere in Platonic dialogues).

Gregory Vlastos, in his famous paper of $1954$ \cite{vlastos1}, pointed out that the TMA  is formally a \textit{non sequitur}
and he investigated the suppressed premises of the argument. There, he proposed the NI and SP axioms, asserting though that
the only explicit premise of the TMA is the \textsl{One-over-Many} (OM) assumption.

\begin{quote}
 ($OM_v$) If a number of things are all $F$, there must be a single Form $F-ness$, in virtue of which
 we apprehend [them] as all $F$.
\end{quote}

Vlastos, after the criticism of his first article (especially by Sellars \cite{sellars}) and a long discussion,
proposed in his second seminal article \cite{vlastos2} a
revised version of the premise-set for the TMA\footnote{Sellars observed that Vlastos \cite{vlastos1}
in stating original NI and SP axioms had used the expression `$F$-ness' as if it represented a proper name of a Form.
Looked at in this way, SP and NI are defective, in that they contain free occurrences of the representative variable
`$F$-ness'. According to Sellars  the defect can be remedied with the aid of quantifiers. Thus, he proposed instead that `$F$-ness'
be taken to represent a quantifiable variable. This simple syntactic maneuver removed the remaining inconsistency.
The TMA's premises as Sellars formulated them in \cite{sellars} are:\\
(OMs) If a number of things are all $F$, it follows that there is an $F$-ness in virtue of which they are all $F$.\\
(SPs) All $F$-nesses are $F$.\\
(NIs) If $x$ is $F$, then $x$ is not identical with \emph{any} of the $F$-nesses by virtue of which it is $F$.
(Actually Sellars' original NI axiom is: If $x$ is $F$, then $x$ is not identical with the $F$-ness by virtue of which it is $F$.
For a criticism on Sellars' original  version see Cohen and Keyt \cite{cohen2} at note $9$). }.

\begin{quote}
($OM_{v2}$) If any of the set of things share a given character, then there exists a unique Form
corresponding to that character; and each of these things has that character by participating in that Form.
\smallskip

\noindent($SP_{v2}$) The Form corresponding to a given character itself has that character (\textsl{Self-Predication}).
\smallskip

\noindent($NI_{v2}$) If anything has a given character by participating in a Form, it is not identical
with that Form\ (\textsl{Non-Identity}).
\end{quote}

There was (and still is) a long discussion on the different ways to express  \textsl{Self-Predication} (SP) and \textsl{Non-Identity}
(NI) premises-axioms. Here we shall adopt  Cohen's \cite{cohen1}, \cite{cohen3} version of them, and we also allow for a `thing' to
be \emph{either a sensible particular or a Form}. On this revised version SP and NI are compatible,
but the three axioms together are not (see Cohen in \cite{cohen1}, \cite{cohen3}). Furthermore, without the
\textsl{One-over-Many} assumption\footnote{Actually we adopt the \textsl{Accurately One over Many} (AOM) assumption,
since we allow for `things' both sensible particulars and Forms. For details see G. Fine \cite{fine} ch. $14$.\label{fnt:aom}}
the theory becomes incomplete. According to many authors, and in my point of view, Vlastos  was
mistaken in supposing that any version of OM with uniqueness quantifier would reintroduce the inconsistency \footnote{This is due since
Vlastos had been working in first-order logic with quantifiers ranging over particulars and Forms.
For more information see also Cohen and Keyt \cite{cohen2} p.$8$ and note $18$.}.

To remedy this problem Cohen in \cite{cohen1}, motivated by Sellars \cite{sellars}, replaced OM with
\\ \textsl{One-Immediately-over-Many} axiom (\textbf{IOM}) and employed more sophisticated logical machinery.
Cohen's reconstruction required quantifying over sets as well (without essentially refuting Vlastos' \cite{vlastos2}
$SP_{v2}$ and $NI_{v2}$ new versions of these axioms). In this new setting Cohen \cite{cohen1} demonstrated
the consistency of the TMA's premises.

Although I do not refute Vlastos' \cite{vlastos1}, pp.$439-440$ textual evidence that the TMA assumes only one Form
$F$ per predicate-character, in this article we follow  Geach \cite{geach}, Cohen \cite{cohen1},
Cohen and Keyt \cite{cohen2} and others in allowing for many Forms per predicate-character in the following context.

We adopt the fact that different things may belong to different levels
(see also Cohen \cite{cohen1} p.$468$, note $31$ and Section $3$ of my article) and based on this,
\emph{from one hand, we  allow for many Forms per predicate, but from the other hand we  argue that these Forms create an increasing sequence
that converges to a unique Form (per predicate)}.

We analyze and defend this thesis, providing the mathematical framework in which this infinite sequence of $F$ Forms is
constructed and it converges-terminates in a unique $F$-Form (per predicate). Furthermore, we address possible problems that seem to arise
from this approach. In particular, within  this framework, we show that that the logic of the TMA is consistent with Plato's intention
for the existence of a unique Form per predicate, that is the uniqueness thesis (see note \ref{fnt:unique}).
Plato for defending the uniqueness thesis introduced in \textit{Rep}. $597c-d$ the famous Three-Bed-Argument, TBA.
An analysis of the logic of the TBA (see Cohen \cite{cohen1}, part VII) definitely calls for the existence of
\emph{at most one} Form per predicate; even though Plato meant TBA to defend \emph{uniqueness} than
`\emph{at most uniqueness}' thesis\footnote{For an elaborate treatment of the uniqueness thesis consult
G. Fine \cite{fine} pp $117, 189-190, 231$.
For its relation to TBA see \cite{fine} pp. $235-238$, and to TMA see \cite{fine} pp. $204, 208-211$.
Also W.K.C. Guthrie \cite{guthrie}, p.$552$, provides another  explanation for
Plato's need to support the uniqueness thesis.\label{fnt:1} \textit{ad.loc.}}.
We will also discuss this in Section $3$B.

Our mathematical interpretation is based on Plato's comprehension of mathematics and particularly
the concepts of \textsl{apeiron, peras, limit, density} and \textsl{convergence} as well as his approach on
\textsl{incommensurable magnitudes}\footnote{ I am not going to accurately translate the terms peras (\textgreek{p\ea ras}) and
apeiron (\textgreek{\asa peiron}) since any translation is an interpretation. My purpose is to understand the meaning
and the characteristics of these terms by the analysis provided mainly in Section $4$.
In general,  \textsl{apeiron} should be understood as \textsl{infinite} or \textsl{unbounded},
related to \textsl{incommensurability}, and \textsl{peras} as rather the \textsl{definite} or \textsl{discrete},
related to  \textsl{commensurability} of magnitudes \label{fnt:peras}.}.
Moreover, our approach is consistent with Plato's philosophy and his Dialectic theory. In supporting our thesis
we provide and  properly analyze concrete textual evidence from Plato's dialogues, Plato's commentators
such as neoplatonists philosophers Plotinus and Proclus, as well as the bibliography developed in this area.

We follow Cohen \cite{cohen1} without planning to raise any objections to $NI_{v2}$ and $SP_{v2}$ axioms.
Essentially these axioms were accepted also
by Cohen\cite{cohen1}, Cohen and Keyt \cite{cohen2} and others. In the sequel we shall also see that a new form of OM axiom
(based on a different construction of TMA), namely the IOM axiom, entails both $NI_{v2}$ and $SP_{v2}$ axioms.
We note that it is not the purpose of this paper to discuss any issues regarding different interpretations of the SP axioms
such as \textit{Broad Self Predication}, \textit{Narrow Self Predication}, or \textit{Pauline Predications}\footnote{For an elaborate presentation, analysis and thorough discussion of the different versions and problems regarding predication
we refer to G. Fine \cite{fine}, ch. $4$ and notes $76$ and $80$. For \textit{Pauline Predications} see also Vlastos \cite{vlastos3}.}.

It is rather worth noting that the \textit{Self-Predication} and \textit{Self-Participation} are completely different concepts.
That is, \textit{Self-Predication} (\emph {SP}) tell us that: \emph{$F-$ness is $F$}, or in other words, is \emph{predicated} as being $F$.

As for NI axiom we  interpret it in the sense that:

\emph{\textsl{`The Form by virtue of which a set of things are all (predicated as) $F$ is not itself a member of that set. Equivalently,
nothing is $F$ by virtue of participating in itself.'}}

Henceforth, the NI axiom, if we adopt that participating in $F-$ness is supposed to explain being $F$,
tells us that we can not explain a thing being $F$ by appealing to this very thing. In this sense NI should be more accurately
phrased as Non-Self-Explanation (NSE) axiom.

We state the versions of SP and NI axioms (slightly modified from $NI_{v2}$ and $SP_{v2}$) that we adopt.

\vspace{0.15cm}

 \textbf{Definition 1}

 \emph{(SP)}: The Form by virtue of which things are (and are judged to be) $F$ is itself $F$. \\
 That is: $\emph{F-ness\,}_i$ \emph{is} \emph{$F$}, $i \in \mathbb{N}$.

 \emph{(NI)}: The Form by virtue of which a set of things are all $F$ is not itself a member of that set.
  Equivalently, nothing is $F$ by virtue of participating in itself, or nothing is explained being $F$ by appealing to this very $F$.\\
  That is: $\emph{F-ness\,}_i$ \emph{does not participate in} $\emph{F-ness\,}_i$, $i \in \mathbb{N}$.

\vspace{0.15cm}

In the sequel, for constructing the TMA we adopt Cohen's approach as developed in \cite{cohen1}.
Our definitions are given in terms of a single undefined relational predicate, `participates in'. The schematic letter `$F$',
which shall serve as a dummy predicate, will play the role that `large' plays as a sample predicate in the TMA.
We note that almost universally it is assumed that Plato intended the TMA to hold for any predicate for which there is a Form;
hence the letter `$F$' is typically used to express this generality.

\vspace{0.15cm}

\textbf{Definitions $\textbf{2}$}:

(D1) By an $F$-object (object for short), we mean any $F$-thing; anything that is, whether a particular-sensible thing
(\textgreek{a\is sjht\og}) or a Form (\textgreek{\is d\ea a}), of which `$F$' can be explained.

(D2) An $F$-particular (hereafter `particular' for short) is an object in which nothing can participates in. That is a
 \textit{sensible thing} (\textgreek{a\is sjht\og}) in Plato's terminology.

(D3) A \emph{Form }is an object that is not a particular.

\vspace{0.15cm}

Cohen's analysis of the TMA clearly aims to exploit an analogy of the TMA with number theory.
It suggests that Plato's infinite regress of Forms (as in the TMA) is analogous to the generation of the infinite sequence
of  Natural Numbers ($\mathbb{N}$) according to  Peano axiomatic construction. It is crucial to state that
Peano's axiomatic foundation of $\mathbb{N}$ (using the concept of successors) is compatible with Plato's infinite sequence of Forms
as well as with the theory of Numbers as developed in Plato's academy\footnote{For example, Plato considers that number $2$
(as a Form) is not the result (or identical) of summation $1+1$, neither the result of the division of a magnitude
in two parts (\textit{Phaedo }$101c$); it is should be considered as the \emph{successor} of $1$.
See also Cohen and Kyet \cite{cohen2}, Fowler \cite{fowler} and Taylor \cite{taylor} ch. $20$.}.

Cohen's TMA construction\footnote{For the development of Cohen's construction see \cite{cohen1} pp.$461-467$ \label{fnt:constr}.}
is in accordance and have counterparts with the standard Von Neumann set-theoretic construction, that corresponds
to the definition of $\mathbb{N}$ via  Peano's Postulates.
The only premise that the TMA seems that does not have an immediate counterpart among the Peano's Postulates is OM. Since OM's role
is to generate a new Form at each stage of the regress, its number-theoretic counterpart could be the \emph{successor} function which generates
the members of the infinite sequence of $\mathbb{N}$. Thus, should OM be appropriately modified, it could be
considered as operating for deriving the principle of Mathematical Induction. Specifically, if OM is to have a uniqueness quantifier,
then it will need to be based on something stronger than Plato's \textit{over} relation. Cohen  \cite{cohen1}
and Cohen and Keyt \cite{cohen2}, understanding this and using Plato's over relation in the new version (IOM) of OM axiom, adequately
define and justify the \textit{immediately-over} function that corresponds to the successor function:

\vspace{0.15cm}

\textbf{\textit{Immediately-Over Function}}

$y$ is immediately over $x \;\; =_{df} \;\; $ $y$ is over $x$, and there is no $z$ such that $y$ is over $z$ and $z$ is over $x$.

\vspace{0.15cm}

It is worth mentioning that the \textit{immediately-over} function is a function \emph{only} with respect to a single sequence of
$F$-Forms. That is, one $F$-Form is `immediately over' another, if no third $F$-Form intervenes between
the two\footnote{The \textit{immediately over} definition can  also be defined in terms of the \textit{over} relation,
the notion of the \textit{level} of an object and the definition of the \textsl{maximal set}, as in Cohen \cite{cohen1}.
Here we use an equivalent one as stated in  R. Patterson in \cite{patter}, p. $54$.
The equivalence of the two definitions is also proved by by R. Patterson in \cite{patter}.}.
Cohen's \textit{One-Immediately-Over-Many} axiom, which entails \emph{OM}-axiom, guarantees that every Form,
in a certain sequence of Forms, has a \emph{unique successor}.

\vspace{0.15cm}

\textbf{(\textit{IOM-axiom)}} \emph{For any set of $F's$, there is exactly one Form immediately over that set}.

\vspace{0.15cm}

This axiom blocks self-participation, since it entails that Forms do not belong to the sets they are over; thus, the axiom NI is built
in IOM-axiom.  In addition, SP (as in Definition $1$) is presupposed as well, because the values of the variables in the definition of \textit{immediately-over} function have been restricted to objects that are $F$.

Specifically, Cohen's construction essentially forms an increasing  sequence $\{F-ness_i\}$, $i \in \mathbb{N}$, of $F$-objects.
There, every object $F-ness_i$, in short $F_i$, defines a level $i$\footnote{See note \ref{fnt:constr}.}.
His construction is in accordance with Von Neumann set-theoretic one\footnote{Recall that in the Von Neumann construction
each member of $\mathbb{N}$ is a member of all its `descendants'.
This fact and the Peano's postulate, stating that `\textsl{no two Natural Numbers have the same successor}',
entail that no member of $\mathbb{N}$ is its own successor. The Form-theoretic analogue of this is that no Form in the sequence
participates in itself.}. Henceforth, the infinite sequence $\{F_i\}_{i=0}^{\infty}$ becomes:
\vspace{-0.2cm}
\begin{align*}
F_0&\hookrightarrow P\\
F_1&\hookrightarrow P\cup \{F_0\}\\
F_2&\hookrightarrow P\cup\{F_0,F_1\}\\
F_3&\hookrightarrow P\cup\{F_0,F_1,F_2\}\\
etc.
\end{align*}

In the above, $P$ denotes the set of $F$ particulars and $\hookrightarrow$ denotes a one to one relation pairing a
Form with the set of its participants\footnote{The symbol $\hookrightarrow$ is used instead of the symbol \,`$=$' \, to avoid
any conceptual identification with the strict mathematical notion of equality or identity.}.
Now using the symbol $\in$ to represent participation rather than set-membership, we obtain the infinite increasing sequence:
\vspace{-0.15cm}
\begin{equation*}\tag{$\star$}
P\in F_0 \in F_1 \in F_2 \in F_3 \in...\in F_k \in... .
\end{equation*}

In our thesis we will deviate (somewhat) from the strict $\mathbb{N}$ analogy.
Namely, \emph{we adopt the above infinite sequence by equipping it with a certain topology}. Under this new perspective
we establish our main result as given in the following claim. The proof of the claim is mainly presented in Section $3$.
\smallskip

\textbf{Main Claim}:

(a) \textit{For every predicate $F$ there exists a countable increasing sequence of $F$-objects, \\
\indent $\{F_i\}_{i=0}^\infty$,  as in} \, ($\star$).
\smallskip

(b) \textit{The infinite  countable increasing sequence $\{F_i\}_{i=0}^\infty$ from part (a)  converges to \\
\indent a unique $F$-object (as a limit) under a certain topology. \\
\indent That is $\displaystyle{\underset{i \to +\infty} \lim F_i=F}$}.
\smallskip

It is worth noting that the above $F$-Form, as a limit of the infinite sequence, could be identified with \emph{the}
unique $F$-Form, the $F$-itself for a certain predicate\footnote{In  Plato's language is refereed also as `F-itself',
 `\textgreek{kaj>{} a\r ut\og}' $F$ (see \textit{Phaedo} $74b3-4$, $100b5-7$ \textit{et.al}).}. Henceforth, it is in accordance
with Plato's intention for supporting the uniqueness thesis (see the previous discussion
and note \ref{fnt:1} \textit{ad.loc.}). Furthermore, since the \emph{least upper bound} is \emph{attained} (as being a mathematical limit),
part (b) of our Claim overrules A. Codes' thesis for the role of TMA as raising an objection to the theory of judgement
(see Section $1$ \textit{ad.loc.}).

The convergence to this limit Form should be understood in the concept of Plato's academy mathematics and it must be comprehended
in a mathematical framework. We note that Plato had a good knowledge of Eudoxus's method of exhaustion for approximating lengths and areas
(see Taylor \cite{taylor}, ch.$20$, `Forms and Numbers' and Anapolitanos \cite{anapol}).
Plato had also used the technique of anthuphairesis (\textgreek{\as njufa\ia resis}) throughout his
work\footnote{ For the technique of anthuphairesis consult
Fowler \cite{fowler}, pp.322-328,  Anapolitanos \cite{anapol} and primarily Negrepontis \cite{negrep}.
These matters are also studied for the needs of our article in Section $4$ \textit{ad.loc.}.}.
In Section $4$ (based primarily on \textit{Philebus}) we shall provide adequate evidence that Plato apprehended
the two different notions of \emph{infinite}, namely the \emph{countable} and \emph{uncountable} ones,
both \emph{in relation to commensurable and incommensurable magnitudes} and the concepts of \emph{density} and the \emph{continuum}
(see also Karasmanis \cite{karasm}). In this sense, it is highly probable that Plato  had also a good grasp on irrational numbers.

\section{The Defense of the Main Claim and Related Issues}
In this section I  defend and eventually provide an adequate proof of my thesis as it was stated in the \textbf{Main Claim}.

I am convinced that Cohen's  TMA version, as presented in the previous section, is closer to Plato's
apprehension of the theory of Forms and to his inclination to think that while the \textit{One-over-Many} axiom
yields exactly one Form for the set under consideration at each step, that principle is consistent with there being more
than one Forms over the set with which we start.
More precisely, over the set of $F-$things just one Form appears or comes into view, even though it turns out
they will  appear more in the process (by repetitive applications of the TMA on the new sets).
Our analysis will set the mathematical framework of the above construction, providing essentially a proof of the Main Claim.
In addition, within our framework, we address, analyze and finally resolve some problems regarding the logic of the
TMA and its relation to Plato's theory of Forms.

It is important to emphasize that if one considers the first time a set $S_{0}$ of $F-$things, \emph{exactly one} $F$-Form $F_0$
shall appear immediately over that set, hence there is exactly one $F$-form.
In the second step, where the TMA is applied on $S_0 \cup F_0$ (the set appeared in the first step and the Form $F_0$)
again one and only one $F$-Form $F_1$ shall appear immediately over this new second set. This process continues up to infinitum,
thus creating an infinite sequence of $F$-objects; the sequence of Forms $\{F_i\}_{i=0}^\infty$.

Next I present the crucial arguments that defend the Main Claim. In subsection \textbf{A} I prove part (a) of the Main Claim
and in subsection \textbf{B} the part (b) of it.
\vspace{0.15cm}

\hspace{1cm}\textbf{A.} At this point we have to clarify what is involved in the claim that a Form can be `over' its participants.
It is clear that Plato thought of Forms as being on a higher ontological level than the sensible particulars participated in them.
Towards this direction there is a strong textual evidence in Plato's work:
\textit{cf., e.g., Rep 515d, 477ff.; Phdo. 74a, 78d ff.; Tim. 28a, 49e; Symp. 210a-212b., et passim}.

As we have already seen, the TMA seems to extend this notion by assuming, in general, that
\emph{a Form is on a higher level than its participants, either these are sensible particulars or Forms}.
That is, each \emph{new} Form that appears in each application of the TMA, see ($\star$),
is in a higher `level'\footnote{We have to stress that these `levels' do not denote different degrees of existence,
since this should be incompatible with Plato's ontological dualism theory.}, within the \textsl{Platonic Realm of Forms},
than its predecessor and hence its participants.
Mathematically, according to ($\star$),  it is formed an increasing infinite sequence of $F$-objects, $\{F_i\}_{i=0}^\infty$.

I defend the above in the grounds that it provides a precise formulation of the logical structure implicit in Plato's arguments.
The key is hiding in how Plato interpreters the \textit{One-over-Many} principle.
There is no indication that
Plato himself ever tried to restrict the \textit{One-over-Many} principle in the way to generate one Form
for each predicate but no more than one. We can support this thesis by recalling \textit{Phaedo's} doctrine
of the homonymy of Forms and their particulars as well as the interpretation of
the famous formulation in \textit{Rep.} $596a6-a7$:

\begin{quote}
{\small\textsl{`We are in the habit of assuming one Form for each set of many things to which we assign the same name}'.}
\end{quote}

According to $596a$, if for a set of many things `we give the same name', then on this set the \textit{One-over-Many} principle
could be applied.  Based on this and in relation to the TMA, Cohen in \cite{cohen1} p.$474$, argues, and I agree with him,
that  it seems inevitable that Plato would ultimately include Forms in sets to which \textit{One-over-Many} principle
is applicable\footnote{See note \ref{fnt:aom}.}.
Here we must also note that the wording of \textit{Republic} $596a$ does not commit Plato to the existence of a Form to
every predicate\footnote{For example, in \textit{Politicus} $262a-e$ Plato denies that there are forms corresponding to every general term.
To know what Forms there are, we need to known not only what words a language contains, but what the genuine divisions in nature are.}.
For Plato may well use the word `name' (\textgreek{\osa noma}) not for every predicate, but for every name
that denotes a property or, as we might say, for every \textit{property-name} in the sense of \emph{true-correct names} as expressed
in G. Fine \cite{fine}, p.$304-305$, notes $44$ and $46$\footnote{This argument is presented clearly in \textit{Cratylus} $386$\textit{ff}
where he seems to use `name' in a more restricted way. According to \textit{Crat.}, `\textmd{n}' counts as a true-correct name only
if it denotes a real property or kind and reveals the outlines of the essence (\textgreek{o\us s\ia a}).
G. Fine in \cite{fine}, pp.$112-113, 304-305, 315$ provides an extensive analysis and a thorough exposition of these notions,
and on what Plato meant by the word `name' (\textgreek{\osa noma}).}.

Additionally, it is worth mentioning that Plato's dialectic, as presented in \textit{Rep}. $509d-511e$, as well as
in $532d-535a$\footnote{ For an analysis of these passages from \textit{Republic} and their relation to Plato's dialectic consult Karasmanis \cite{karasm2}, \cite{karasm3}, Benson \cite{benson}, J. Annas
\cite{annas} ch $10, 11$, as well as Robinson \cite{robinson} ch. $6,7,10$.\label{fnt:dial}},
is in the general line of thought of part (a) of the Main Claim. 
In particular, Plato in the analogy of \textsf{Simile Divided Line}  advances the hierarchical progressing model of levels
that leads to the \textsf{anupotheton arx\={e}n}, claiming that the whole procedure is done in the realm
of \textsl{no\={e}sis} (\textgreek{n\oa hsis}) and precisely in the section of
\textsl{epist\={e}m\={e}} (\textgreek{\es pist\ha mh}), via the exercise of dialectic method.

Analogous arguments are advanced in the seminal work of Proclus' \textit{Commentary on Plato's Parmenides} (translated by  G. Morrow)
and particularly in  $879. 15-28$ and more emphatically in $881. 23-33$\footnote{We quote: \textit{`...And from there in turn he will
be chasing after unities of unity, and his problems will extend to infinity, until, coming up against the very boundaries of intellect,
he will behold in them the distinctive creation of the Forms, in the self-created, the supremely simple, the eternal...'}.}.

It is worth mentioning that in \textit{Symposium} the concept of  `\textsl{Rungs}' of `Diotima's ladder'
and Plato's  analysis on this support the model of increasing sequence of Forms. There, the procedure is presented in a detailed and vivid manner
where the limit Form $F$ is the one of `Beauty'. Of course there the corresponding convergent sequence $\{F_i\}$ consists
of $F$-objects that are not Forms. This fact is not of a major importance since in this paper the references to \textit{Symp.}
are primarily aiming in clarifying and analyzing the procedure itself as well as  the concepts of density and convergence.

Particularly, Plato in \textit{Symp.} $210e$, among others, states
\textit{`...passing from view to view of beautiful things, in the right and regular ascent,..'},
noting also that the ascension to the final Rung, corresponding to the ultimate Form (of `Beauty') itself, has to be done in a
\textit{`correct and orderly succession'} (\textgreek{\es fex\hc s {}\<orj\wc s t\ag{} kal\aa}).
This is even more clear in \textit{Symp.} $211b-c$ where the nature of the ascending procedure to the true-Form of Beauty is analyzed.
From this passage we hold on the phrase `\textgreek{ \wra sper \es panabasmo\ic s {}qr\wa menon}', `as on the rungs of a ladder',
and the use of the word `\textgreek{\es pani\wa n}', `that ascends'. According to Plato this describes
the `right approach' (\textgreek{\os rj\wc s}) for `almost being able to lay hold of the final true $F$-Form'
(`\textgreek{sqed\og n{} \asa n{} ti{} \ara ptoito{} to\uc{} t\ea lous}'), which constitutes also the ultimate goal
and the conclusion of the whole procedure. There, each \textit{rung of the ladder} defines a level (in an analogy to ($\star$)),
where each level is higher and contains the preceding ones. According to Vlastos \cite{vlastos4} the whole procedure moves
`closer step by step to the Beauty itself'.
We have to state that nothing prevents us from assuming that the same model holds analogously for all Forms (predicates)
and is not restricted only to the Form of `Beauty'.

The ascending procedure described in \textit{Symp.} could be considered as analogous to the abstract one presented in \textit{Rep}.$511b$,
where we hold the phrase `\textgreek{o\irc on {}\es pibib\aa seis}'. This ascending procedure it is a fundamental one within
the Platonic dialectic (see also note \ref{fnt:dial}.).

It is convincing that the previous analysis strongly supports the
concept of \emph{degrees of hierarchy} among the plurality of Forms $F_i$, $i \in \mathbb{N}$, that appear
by repetitive applications of TMA. Thus, the sequence $\{F_i\}_{i=0}^\infty$ (in ($\star$)) is justified as increasing.

\vspace{0.15cm}

\hspace{1cm}\textbf{B.}  I  establish the second part of the main Claim. Thus, I show that  the infinite regress
of $F$-Forms  $\{F_i\}_{i=0}^\infty$ (for a particular predicate) \emph{converges} to the unique (terminating) $F$-Form,
the so called  $F$-\textsl{itself} (\textgreek{kaj>{} a\r ut\og} `$F$') for the predicate in concern.
Furthermore, I  address the various questions that arise regarding the topological framework of this convergence, as well as the nature of the limit-Form  $F$ and the way it should be understood in relation to the various $F_i$, $i\in \mathbb{N}$, of the sequence in concern.

I will argue that the TMA model adopted in this article, based in Cohen's analysis \cite{cohen1},
is also compatible with Plato's intention that \textit{One-over-Many} principle (in its \textit{One-Immediately-over-Many} version)
yields to a uniqueness thesis. In this framework we shall defend that the logics of the TMA and the TBA are not inconsistent, but consistent
and rather mutually complementary. Thus, the \emph{uniquness} thesis, that Plato intended to support by introducing the TBA,
should be apprehended in the context that the unique $F$-Form is the limit of the convergent sequence $\{F_i\}_{i=0}^\infty$; that is,
$\displaystyle{\underset{i \to +\infty} \lim F_i=F}$\footnote{This $F$ is the unique Form, as defined in \textit{Phaedo},
the `F itself ' ($74a11-12$), the `F without qualification' ($74b$\textit{ff}), the ` `F' that it can never seem non-`F''($74c1-3$),
see also note \ref{fnt:itself}. It is \emph{the} Form as the final stage of the ascending procedure described above.
In \textit{Symp.} $211c$ this $F$ is further understood as the unchangeable end, the goal, the conclusion of the
ascending procedure, `\textgreek{a\us t\og{} teleut\wc n{} \org{} \esa sti}', tangent to the very essence of $F$-ness.
Similar terminology and way of apprehending this `$F$-itself' is encountered in many Platonic dialogues, such as
\textit{Phaedo} \textit{et.al}. For example in \textit{Phaedo} $101e$ it is described as the termination of the ascending procedure
to the \emph{one} Form which is `adequate', `\textgreek{\era ws{} \es p\ia{} ti{} \ir kan\og n{} \esa ljois}'.
The ascending procedure is also developed elaborately in an abstract manner in \textit{Rep}.$509e-511d$
and shall be discussed in the sequel. See also Karasmanis \cite{karasm2} for the analysis of hypothetical method in this passage.}.

We have to inform the reader that  a different approach on this subject, claiming that the logics of the TMA and the TBA
are rather inconsistent, is presented in G. Fine \cite{fine}, $pp. 235-238$. In defending her thesis, G. Fine
also argues that, in this particular case, Plato consistently rejects the NI axiom.

We stress that Cohen, while examining the consistency and logic of the IOM axiom, discovered that in order for his IOM axiom to be consistent
the set theory his formalization presupposes cannot include the \emph{Principle of Abstraction}\footnote{It is out of the scope
of this article to argue about the intrinsic of the Principle of Abstraction \textit{per se}. For the validity and its difficulties
see Quine \cite{quine},  pp. $134-136, 249, 300$.}:

\vspace{-5mm}
\begin{equation*}
(\exists \alpha)\, (\forall x) \left(x\in \alpha \leftrightarrow Fx\right).
\end{equation*}

That is, for any predicate $F$, there is a set $\alpha$ consisting of \emph{all and only} objects to which
that predicate applies\footnote{This principle should be considered in relation to  \textsl{the axiom schema of comprehension}
of Zermelo-Fraenkel set theory. The interested reader should also study the nature of this axiom schema and the possible
problems arising by an improper use of it. For details and an analysis we refer to Cori and Lascar \cite{cori} pp.$112-113$.}.

Cohen's  problem  was the existence of a universal set of $F$'s; that is a highest level set containing all the $F$-objects
in all (lower) levels. For, if there were such a set (the universal set of $F$'s), the increasing sequence in ($\star$), that corresponds
to the Von Neumann set theoretic construction of $\mathbb{N}$, it would contain a maximal element (and hence should not be an infinite one).
Of course this cannot happen, since it contradicts the IOM axiom and the fact that the set $\mathbb{N}$ does not have a maximal element.

Cohen addressed this problem and he stated that if someone wants to retain the \emph{Principle of Abstraction},
the IOM should have been somewhat altered (see Cohen \cite{cohen1}, note $33$).

Here we follow a different approach, offering a solution to the problem and \emph{retaining the original version of IOM axiom}
\emph{and the Principle of Abstraction}. This is done by claiming that the increasing infinite sequence $\{F_i\}_{i=0}^\infty$ \emph{converges}
within a certain topology; hence  there exists a (mathematical) limit of this sequence. This limit as the attainable
\emph{least upper bound} of the increasing sequence could be considered as the the `highest level set' containing \emph{all} the $F$-objects
in \emph{all} (lower) levels; thus the \emph{Principle of Abstraction} is retained.
We believe that our approach is closer to Plato's theory of Forms and specifically to his intention for accepting
the \textsl{uniqueness thesis}. In the sequel we support this thesis.

Plato in \textit{Rep}. $597c-d$ presents the Three Bed Argument and he applies the \textit{One-over-Many} principle
to  eventually prove the uniqueness of the Form (the `bed' in the particular case).
But as analyzed by Cohen \cite{cohen1} and G. Fine \cite{fine},
the application of the TMA  on the TBA what does really proves is that there is at \textit{most} one Form of `bed'.
It is remarkable that in order to conclude that there is \emph{exactly one} Form of `bed' we must show that the sequence as constructed,
using repeatedly  the TMA, eventually stops. This cannot be done, since the TMA produces infinite sequence.
Thus, the logic of the TMA together with that of the TBA lead to the conclusion that there are \emph{none} Forms.
Of course this  was not the intention of Plato since it contradicts his Ontological theory.

Here we provide adequate evidence for establishing the proof of part (b) of the main Claim. This also entails that that
the \emph{Principle of Abstraction} and the IOM are retained,  and the existence of the infinite sequence
is proved to be genuine without leading to any contradiction. In addition, within this context,
A. Codes' thesis that the TMA raises an objection to the theory of judgement (see Section $1$ \textit{ad.loc}) is overruled.
For doing this we argue that the $F$-Form, the `F-itself', that appears as the limit-terminating point of the above increasing sequence
must not be committed to the \textit{Non-Identity} (NI) or \textit{Non-Self-Explanatory} (NSE) axiom .

Indeed, form one hand  this $F$-Form clearly satisfies Self-Predication axiom (since it is predicated as being $F$).
But from the other hand, it is the limit of the \emph{increasing} sequence of $\{F_i\}_{i=0}^\infty$ and hence there
are no further $F$-Forms beyond this particular $F$-Form, in contrast with the rest $F_i$, $i\in \mathbb{N}$,  $F$-Forms
in the sequence ($\star$). These arguments lead to the conclusion that the limit $F$-Form
\emph{should be comprehended as self explained}. Hence it could not satisfy the Non-Identity (or NSE) axiom.

In the sequel we present how Plato understands and explains the above thesis. For doing this
we study the framework within which he comprehends, in my point of view, the convergence of $\{F\}_{i=0}^\infty$
to the unique  $F$ (for each predicate).

Of course in order to talk about convergence and limits we must assume a certain topology.
This topology should be the one closer to Plato's understanding of mathematical concepts. Plato, as we mentioned in Section $2$
and as we  shall see in Section $4$, had a good grasp on the Eudoxus' exhaustion method and the technique of anthuphairesis,
as well as an apprehension of the concepts of \textsl{peras, apeiron, density and continuum} and their relation
to \textsl{commensurable} and \textsl{incommensurable} magnitudes.

The topology we consider, and hence the convergence we are refereing to, is also in accordance with  Plato's dialectic process
(more specifically as presented and analyzed in \textit{Rep.} $509c-511e$ and $532d-535a$) as well as his apprehension of the
fundamental concept of the \textsf{anupotheton arx\={e}n} (\textgreek{\as nup\oa jeton{} \as rq\hg n}) in \textit{Republic},
but applied here for each particular predicate. In \textit{Rep.} $510b$ Plato states clearly that the highest rung of the ladder
is not reached until the entire domain of  epist\={e}m\={e} has been exhausted via the dialectic process. This principle, established
by an exhaustive scrutiny, should be understood as higher than all premises-hypotheses, `\textgreek{\ur poj\ea seis}'.
It is higher, in the sense that contrary to them it has an axiomatic status (playing the role of a system of axioms),
it is non-hypothetical, it is situated in the highest point of the intelligible world (in \textit{Republic's}
\textsl{Simile Divided Line}) and it does not require derivation\footnote{As argued above even though it is not committed to NI axiom,
it is predicated as being  $F$ and thus it satisfies the SP axiom.} (see also Karasmanis \cite{karasm2}, \cite{karasm3}
and Benson \cite{benson}, p.$190$).

We emphasize that mathematically the \textsf{anupotheton arx\={e}n} is the ultimate-final unique Form apprehended
as the mathematical limit of the infinite increasing sequence under the presupposed  topology.
It must be noted that this Form is comprehended not as a transcendental ontological mystery
but in the mathematical sense of the \emph{least upper bound} of the increasing sequence of Forms (see ($\star$)
and part (a) of the main Claim) that it is eventually attained; hence it becomes a limit Form.
Apart from \textit{Republic}, it should be apprehended as  the \textbf{one} (`\textgreek{monoeid\eg s}' in \textsl{Phaedo's} language)
that should be parallelized with the highest-terminating rung of \textsl{Diotima's} ladder, which is tangent to the
very essence of $F$-ness\footnote{ see \textit{Symp.} $211b-c$ where that $F$ is the Form of `Beauty'.}.

Within  this framework, this limit-Form could also be conceived as not committed to the Non-Identity (NI) axiom.
This is due for being the $F$-Form, the \textsf{anupotheton arx\={e}n} that in addition to the above gives an account to all the
lower level $F_i$, $i \in \mathbb{N}$, $F$-Forms, but itself does not require derivation.
In this sense, there do not exist further $F$-Forms in higher order level(s) to be depended on,
or any Forms for providing an explanation to it\footnote{Strictly, this should be better considered as modal, stating:
`that there are not required further $F$-Forms in higher order level to be depended on,
or needed for providing an explanation to it.' But no harm is done here by simplifying it as existential,
since the TMA shows existence of $F$-things. Furthermore, in a mathematical-logic language, this unique Form is not a derivation,
or a theorem, but it has the status of an axiom, or of a system of axioms,  and thus it does not require a proof.
Moreover, due to its status, all the information of the system can be retrieved from it \label{fnt:axiom}.}.
Thus, the Form in discussion, though is predicated as $F$, should be regarded as self-explanatory.
In an abstract and an ontological level this Final stage (this particular limit of the sequence) should be considered as
the unchangeable unique `\emph{Is}'\footnote{In \textit{Symp}. $211a$ Plato characterizes this `\emph{Is}' as:
`\textsl{... ever-existent and neither comes to be nor perishes, neither waxes (growths) nor wanes (declines, decreased)...}'.
In $211b$ is characterized as unchangeable and is affected by nothing. Further, in  $211c$ this `\textsl{Is}' is revealed
at the end of the ascending procedure, characterized as the very essence of the $F$-ness. },
`\textgreek{e\isc nai}', predicated as  $F$ and apprehended as the terminating point-Form
of the increasing infinite sequence-process, the process of \emph{becoming}-gignesthai ( `\textgreek{g\ia gnesjai}').
This infinite process, using Plato's terminology, leads to the \textsl{unconditional}, \textsl{immutable},
\textsl{objective}, \textsl{unchangeable}, \textsl{perfect} and \textsl{unique} `\textit{Is}'.

In Section $4$ we present and analyze, mainly within the context of \textit{Philebus} $23c-27c$, the mathematical framework of this process.
It is the process \textsl{per se} which is compactly phrased at $26d7-9$ as \textsf{genesis eis ousian}-
`\textgreek{g\ea nesis{} e\is s{} o\us s\ia an}'\footnote{For an analysis of this passage, as well as the crucial phrase
`\textgreek{g\ea nesis{} e\is s{} o\us s\ia an}' we refer also to Section $4$.}.
This procedure it can be parallelized to the \emph{dialectic} one, that under the cause of \textsf{no\={e}sis} exhausts
the entire domain of knowledge-\textsf{epist\={e}m\={e}} in order to arrive and terminate to the purest mental state,
the unique $F$-Form. This Form, as stated earlier, is non-hypothetical (`\textgreek{\as nup\oa jeton{} \as rq\hg n}'
using \textit{Republic's} terminology) is not committed to NI (better NSE) and has an axiomatic status.

Summarizing, we conclude that the increasing infinite sequence of $F$-Forms derived
by repetitive applications of the TMA (for each predicate) converges eventually to a (mathematical) limit
(a terminating point) which is unique (as a limit of a sequence).
This is in accordance to Plato's intention in supporting the  \emph{Uniqueness of a Form} per predicate.
Furthermore, our thesis is consistent with the \emph{Principle of Abstraction} in Logic as analyzed earlier in this paper.

At this point we provide further textual evidence supporting the concept of \emph{convergence} (as stated earlier)
strengthening  our arguments in part (b) of the main Claim.

We mention that the notion of the \textit{limit}, as the terminating point of this infinite but convergent sequence,
as well as the whole theory established above, it is also present in some form or as parts of it in many Platonic texts
(of the middle and late period), such as \textit{Symp.}, \textit{Phaedrus}, \textit{Rep.},\textit{ Philebus},
 \textit{Epistle} $7$  \textit{et.al}. Its distinct status, primarily in \textit{Philebus}, shall be analyzed in Section $4$.
In addition, the notion of the \textit{limit} was adopted and studied by Neoplatonists philosophers and Plato's commentators
such as  Proclus  and Plotinus.

More specifically, Plato in \textit{Symp.} $210e-211a$ characterizes the limit Form $F$, in relation to
the termination of the infinite sequence, as being revealed `abruptly, suddenly'-`\textgreek {\es xa\ia fnhs}'.
In addition he states that it exist unconditionally and is `\textit{the perfect thing, the wondrous and beautiful in nature}'-
`\textgreek{jaumast\og n{} t\hg n{} f\ua sin {}kal\oa n}'\footnote{A similar terminology is used in \textit{Phaedrus} $250b$.}.
Furthermore, in Plato's \textit{Epistle} $7.341c$ we encounter related notions that emphasize the concept of the
\textsl{upper bound of the infinite sequence that it is revealed suddenly}-`\textgreek{suz\hc n{} \es xa\ia fnhs}'\footnote{\textit{Epistle}
$7341c$ `...\textsl{but, as a result of continued application to the subject itself and communion therewith,
it is brought to birth in the soul on a sudden, as light...}'.}.
Analogous approach is encountered  also in Plotinus \textit{Enneads} $43.17$.

The statement that the terminating $F$-Form can be achieved as a limit  of the increasing sequence
in the \textsl{most perfect manner} is furthermore emphasized in the process described in \textit{Symp.} $211a-e$.
There, it is characterized (in $211e$) as  `the divine beauty itself, in its unique form'-
`\textgreek{je\ic on{} kal\og n{} monoeid\eg s}'\footnote{In \textit{Republic} $398a$ this perfect from
is characterized as `\textsl{divine and holy}'-`\textgreek{\ir er\og n{} ka\ig{} jaumast\oa n}'.},
and in  \textit{Symp.} $212a$  as the \textit{tangent-contact to the truth}-`\textgreek{to\uc{} \as lhjo\uc s{} \es faptom\ea n\wi}'.
In \textit{Symp.} $211c9$ is called  the `terminating point of the ascending procedure'-
`\textgreek{a\us t\og{} teleut\wc n{} \org{} \esa sti}'\footnote{For an analysis of the passage $201d-212c$
in \textit{Symposium} we recommend Taylor \cite{taylor}, ch $9$, section $8$ and Vlastos \cite{vlastos4}.}.

Similar ideas about the nature of this `limit', the terminating point-Form,
as being the contact approach and intercourse with the truth is evident throughout the Platonic corpus\footnote{This idea is
present and analyzed in many Platonic works and is also frequent in Aristotle and the Neoplatonists.
The process of arriving to this terminating point that is tangential to the true-Form is analogous to the one described
earlier in \textit{Symp.} $210e$\textit{ff}, \textit{Rep.} $490b$, as well as  in \textit{Rep.} $509c-511d$ in the context of
Plato's dialectic theory. (For details see J. Annas \cite{annas} and for a meticulous analysis
of \textit{Rep.} $509c-511d$ see Karasmanis \cite{karasm3} and Benson \cite{benson}).}.

Very emphatically, the Neoplatonist philosopher Proclus in his work \textit{Commentary on Plato's Parmenides}
$881, 23-33$\footnote{We quote:`...\textsl{And from there in turn he will see other more comprehensive unities, and he will
be chasing after unities of unity, and his problems will extend to infinity, until, coming up against the very boundaries of the intellect,
he will behold in them the distinctive creation of Forms, in the self-created, the supremely simple, the eternal...}'.
A similar line of thought is present also in Plotinus, \textit{Enneads} $2.4.15, 15-16$.}
analyzes the concepts of infinite (\textit{apeiron}) regress that arrives  to a terminating
\textsl{mental} point/\textsl{peras}-`\textgreek{noer\og n{} p\ea ras}', via the process of the
intellect/\textsf{no\={e}sis}-`\textgreek{n\oa hsis}'. The term `\textgreek{noer\og n{} p\ea ras}' should be apprehended
as analogous to `\textgreek{\as nup\oa jeton{} \as rq\hg n}' of \textit{Republic}. It is the termination, the limit of the
dialectic process. This process occurs in  the \textgreek{\es pist\ha mh}-\textsf{epist\={e}m\={e}}
(the upper part of \textgreek{n\oa hsis}-no\={e}sis) of the `\textit{Simile Divided Line}'\footnote{For further
clarification and analysis of these terms we recommend  Karasmanis \cite{karasm2}, \cite{karasm3} section III
(\textit{cf}. note $3$ \textit{ad.loc}.) and Benson \cite{benson}.}.
In addition, Proclus in \textit{Commentary of the First Book of Euclid's Elements} is arguing, using
the concepts of \textit{peras} and \textit{apeiron} (in the line of thought of Plato's \textit{Philebus}), in order to establish
the convergence apprehended via the notion of the \emph{limit}.
This approach shall be discussed in some extend in the next section, where we shall compare and cross-examine it with
the one that appears in  Plato's \textit{Philebus}.

In the next section we study, in some extend, the mathematical concepts and the general framework involved and required
for establishing our main Claim. Our analysis is based primarily on Plato's \textit{Philebus }and the key notions of
\textsl{peras}, \textsl{apeiron}, \textsl{density} and \textsl{continuum} in relation
to \textsl{commensurable} and \textsl{incommensurable} magnitudes.

\section{The Mathematical framework of the Argument}
I would like to start by arguing that Plato apprehended the structure of incommensurable magnitudes in a way that these magnitudes
correspond in a unique and well defined manner to the modern concept of the \textsl{Dedekind cut}\footnote{This is more clear in
the dialogues \textit{Theaetetus}, \textit{Philebus}; for details see Taylor \cite{taylor},
ch. $20$, Fowler \cite{fowler}, and Anapolitanos \cite{anapol}.}.
For the  precise definition of the Dedekind cut the reader should consult Rudin \cite{rudin}\footnote{In Rudin \cite{rudin} the
Real numbers are constructed as Dedekind cuts in a unique and well defined manner
which follows  the line of thought of Eudoxus' exhaustion method. It is important to note that square roots of non square numbers,
(see \textit{Theaet.} $147d4-148b4$), as well as the incommensurable magnitudes could be obtained
via the technique of anthuphairesis, `\textgreek{\as njufa\ia resis}'. Antuphairesis entails also the concept of `cut'
(see Taylor \cite{taylor}, ch. $20$) and is presented formally and in great detail in the 10th book of the Elements of Euclid.
See also note \ref{fnt:anth}.}.

Furthermore, it is important to note that Plato  captured the notions of \textsl{density} of magnitudes and
the \textsl{complete continuum} in the sense that they include incommensurable cuts (see Karasmanis\cite{karasm} p.$394$).

It seems that Plato considers incommensurability as an essential feature of magnitudes. This can be most vividly observed
in \textit{Philebus} $23c-27e$, as well as \textit{Theaetetus} $147d4-148b4$, where  Plato's approach to the \emph{continuum} is developed.
In  \textit{Philebus} Plato makes an effort to explore the relation between \textsl{continuum}, \textsl{infinite divisibility} and
\textsl{incommensurability} in contrast with \textsl{commensurable} things that are capable of appearing in ratios and proportions.

We study  the above mathematical concepts and the notions of \textsl{peras} and \textsl{apeiron}\footnote{see note \ref{fnt:peras}.}
based primarily on the  passage  $23c-27c$ of \textit{Philebus}.

We have to note that in $16b-19a$  Plato says that `the things that are ever said to be' (\textgreek{\as e\ig{} leg\oa mena})
are made up of \textsl{one} and \textsl{many}, with \textsl{peras} and \textsl{apeiron} inherent (\textgreek{s\ua mfuton}) in them,
but in the second passage Socrates asserts that \textsl{peras} and \textsl{apeiron}
are two different general classes of things\footnote{We must note that this is indeed a problematic passage where
the term `\textgreek{\as e\ig{} leg\oa mena}' are most probably the ideas-Forms. For a discussion and the various
interpretations of the term \textsl{apeiron} in this passage,
we refer to Karasmanis \cite{karasm} p.$390$, notes $8,9$. Analogous to this approach,  for the so called
\textsf{ont\^{o}s on}-`\textgreek{\osa ntws \osg n}', is advanced by Proclus in \textit{The Elements of Theology}, $89$.}.

In $23$\textit{ff} Plato makes a fourfold division: \textsl{peras} (\textgreek{p\ea ras}), \textsl{apeiron} (\textgreek{\asa peiron}),
\textsl{mixed} (\textgreek{meikt\oa n}) and \textsl{cause} (\textgreek{a\is t\ia a}).
That is how he examines these four categories-classes one by one.

In $23c10-d1$ the class of \textit{mixed}  is assumed as a combination of \textsl{peras} and \textsl{apeiron}.
Karasmanis in \cite{karasm} examines systematically the characteristics of \textsl{peras} and \textsl{apeiron}
and attempts to answer a plethora of questions regarding their nature, their status in the fourfold division and their relation
with the notions of \textsl{commensurability} and \textsl{incommensurability}.

The class of `cause' (\textgreek{a\is t\ia a}) is explained in $23d6-7$ as the cause of the existence of the third class,
that of the \textsl{mixed}. According to Plato this is the cause of the combination (\textgreek{s\ua mmeixis}) of the other two
classes, the \textsl{peras} and the \textsl{apeiron}.

In $23e$ and $25d11-e2$ Socrates states that he should investigate the mechanism (\textgreek{a\is t\ia a}) that the
separated \textsl{peras} and \textsl{apeiron} are mixed together, explaining how are becoming a unity.
Plato in $24c-d$  explains more accurately what he means by `\textsl{apeiron}' and `more and less',
`\textgreek{t\og{} pl\ea on{} ka\ig{} t\og{} \esa latton}'.
He says that wherever the `more and less' are present they exclude any definite quantity, \textit{poson} (\textgreek{pos\oa n}).
This passage says that the essential characteristic of the `more and less' is the absence of any definite quantity,
for the presence of definite quantity and measure (\textgreek{m\ea tron}) in the place where the `more and the less' is present
will abolish the `more and the less'. It seems that according to Plato the notions \emph{`more and less'}
and \emph{`definite quantity'} are mutually excluding.

In $24d$ the notion of \textsl{continuum} is further analyzed. Specifically in $24d4-5$ Plato uses the expression
`\textgreek{proqwre\ic{} g\ag r{} ka\ig{} o\us{} m\ea nei}' translating as `goes on without pause' or
`progressing and never stationary'; henceforth the notion of \textsl{a continuous motion} is advanced.
In the same line of thought  Plato at $31a$ says that:

\begin{quote}
 `\textsl{apeiron in itself does not have and will
never have any precisely marked beginning, middle or end}'.
\end{quote}

Thus, we could assert that what is characterized by the `more and less' should be \emph{continuous} and, therefore,
the characteristics that \emph{Plato attributes to the \textsl{apeiron} point to continuous magnitudes}.
Such an \textsl{apeiron} as continuous is infinitely divisible; this last property is the main characteristic of magnitudes.
The same idea is also clear in  $24e-25a$.
In passages $25b$ and $25e$ he talks about the class of \textsl{peras} where he obviously relates that
to the class of \emph{commensurable} magnitudes (analogously one could relate it to the countable set of rational numbers).
Thus, it seems that Plato suggests that the main characteristic of his \textsl{peras} is \emph{commensurability}.

We find further evidence for this claim in two other passages from the \textit{Philebus}:

(a) `\textsl{That of equal and double, and whatever puts an end to opposites being at odds with each other,
and by the introduction of number that makes them commensurate-summetra (\textgreek{s\ua mmetra}) and harmonious}', $25d11-e2$.

(b) `\textsl{Again, in the case of extremes of cold and heat its advent removes what is far too much and apeiron and produces what is
measured-emmetron (\textgreek{\esa mmetron}) and commensurable (summetron)}', $26a7-9$.

If \textsl{peras} is what makes things \textsl{commensurate}, then \textsl{apeiron} must be the source
of \textsl{incommensurability}. I think that Plato here is using the term commensurate-summetron
in a rather technical, mathematical sense. Now, if we agree that Plato says in  $24e7-25a2$
\footnote{We quote: `\textsl{All things which appear to us to become more or less, or to admit of emphatic
and gentle and excessive and the like, are to be put in the class of the infinite as their unity...}'.}
and in $25a6-7$\footnote{ We quote: `... \textsl{and the things which do not admit
of more and less and the like, but do admit of all that is opposed to them...}'.} that \textsl{apeiron}
admits opposite characteristics to those of \textsl{peras}, then we have to conclude that \emph{incommensurability
is a further and very important characteristic of \textsl{apeiron}}.
It seems then that Plato relates \emph{discontinuity to commensurability} and probably (ex silentio) \emph{continuity to incommensurability}.

I  refer to $26a$ to stress that according to Socrates, the perfection (in the art of music) or harmony and moderation (in the
case of temperature) can be achieved by properly combining-mixing the opposite directions of \textit{apeiron}
(see Karasmanis \cite{karasm}, p.$391-392$) and \emph{finding the limit-peras}\footnote{We quote $26a2-4$:
 ...`\textsl{And in the acute and the grave, the quick and the slow, which are unlimited, the addition of these same elements
 creates a limit-peras and establishes the whole art of music in all its perfection, does it not?}... .}. This limit-peras  within
 this context is restated (and identified) in $26d8-9$
as the emmetron-measured and summetron-commnsurable\footnote{See the notes \ref{fnt:metron} and \ref{fnt:limit}.}.
In $26b$ Plato generalizes the  assertion that the limit leads to perfection and that is
the way to arrive to the realm of Forms (\textgreek{\is d\ea es})\footnote{We quote from $26b$..\textit{`all the beauties
of our world arise by mixture of the apeiron with the peras'}, and he continues to state that \textit{`for  many glorious
beauties of the soul this goddess,... beholding the violence and universal wickedness which prevailed,
since there was no limit-peras \textsf{(see note \ref{fnt:limit})} of pleasures or of indulgence in them,
established law and order, which contains a limit-peras \textsf{(see  note \ref{fnt:limit})}...}'.}.
Now \textit{peras} gives not just any determination of \textit{apeiron} but the right one according to specific case of art.
When Plato speaks about \textit{apeiron} in the framework of mixed and the mechanism
of \textit{mixis}, it seems that he broadens its meaning. In $26d7-9$ Plato states:

\begin{quote}
\small{`...\textsl{And as to the third class, understand that I mean every offspring of these two [peras and apeiron]
\textsl{which comes into being} (\textgreek{g\ea nesis{} e\is s{} o\us s\ia an}), to a stable and immutable essence, as a result
of the limits-measures (\textgreek{m\ea trwn})\footnote{Here the word `\textgreek{m\ea tron}' (measure) should be
comprehended as `due measure or limit'-`right proportion' \label{fnt:metron}.} created by the cooperation of the
the peras\footnote{Here \textsl{peras} should be comprehended rather as `that which limits or has limits'.
See also \textit{Philebus} $30a$ \label{fnt:limit}.}....'}}
\end{quote}

Thus, according to the above passage and to $23d6-7$ (commented earlier) \textsl{peras} is drown onto \textsl{apeiron}
via the class of the \textsl{cause} (\textgreek{a\is t\ia a})
to create the third class, the class of \textsl{mixed} (\textgreek{meikt\og n}) according to a specific process.
Summarizing all of the above we conclude:

\begin{quote}($\maltese$)
The incommensurability of \textit{apeiron} can be approximated and eventually be a limit in a specific manner;
that is, by imposing \textsl{peras}, which makes things \textsl{commnsurate},  on it.
\end{quote}

We could further conclude that Plato considers the class of \textsl{mixed} in relation to the
unchangeable perfect `\emph{Is}'\footnote{See also the discussion in Section $3$ \textit{ad.loc.}.}
(\textgreek{e\isc nai}) which is formed by imposing, with the aid of the class of \textsl{peras},
a limit, a due measure on the class of \textsl{apeiron} via  a specific process. This process that results to the `\emph{Is}',
is characterized in Plato's language ($26d8$) as `\textsf{genesis eis ousian}'-`\textgreek{g\ea nesis{} e\is s{} o\us s\ia an}'.

Plato clearly supports  the above thesis in other parts of \textit{Philebus} such as $27d6-10$:

\begin{quote}
{\small\textsl{`for that class is not formed by mixture of any two things, but of all the things which belong to
the \textsl{apeiron}, bound by the \textsl{peras}; and therefore this victorious life
(\textgreek{\<o{} nikhf\oa ros{} b\ia os}) would rightly (\textgreek{\os rj\wc s}) be considered a part of this class.'}}
\end{quote}

For an in depth comprehension  of this  thesis the translation of the phrase `\textgreek{g\ea nesis{} e\is s{} o\us s\ia an}' in $26d8$
is crucial. It is translated as `coming into being', that entails a continuous process, a creative procedure;
it is a becoming that generates and eventually gives rise to something stable via the proper-correct way in due measure.
In this sense, it leads to the ideal-optimum `victorious life' which belongs to the class of mixed\footnote{Earlier in $27d1-2$
this `victorious life' is characterized as the mixed life of pleasure-h\textsf{\={e}don\={e}} and prudence-\textsf{phron\={e}sis}.
See also  Taylor \cite{taylor}, ch $16$.}. At \textit{Philebus} $55a3$, the \textsf{genesis} as the becoming, the generation, is presented
as the opposite of `destruction'-`\textgreek{fjor\aa}'. The word \textsf{genesis} appears often in Plato. Particularly,
in  \textit{Phaedrus} $245e2-5$ is presented as the source of motion and origin, where the self-motion is comprehended
as the essence, `\textgreek{o\us s\ia a}'. The other term ousia-`\textgreek{o\us s\ia a}' is central in Plato's philosophy.
It primarily denotes the essence, the true substance, the stable and true being, the immutable reality, the
\textsf{einai}-`\textgreek{e\isc nai}'\footnote{In contrast to `\textgreek{o\us s\ia an}' as the `\textgreek{e\isc nai}',
the `\emph{Is}', the  stable being and the immutable reality, Plato uses in \textit{Theaet.}$185c9-10$ the term
`\textsf{m\={e} einai}'- `\textgreek{m\hg{} e\isc nai}'  \label{fnt:2}.}. We have to mention that in \textit{Rep.}$534a2-4$
Plato claims that the two terms-concepts, \textsf{genesis} and \textsf{ousia} have different ontological and epistemological status.
He emphasizes that \textsf{doxa} (\textgreek{d\oa xa}) is concerned and is dealing with \textsf{genesis},
while \textsf{no\={e}sis}-\textgreek{n\oa hsis} is the one that deals with the essence-\textsf{ousia}\footnote{For details about
the status of these terms we refer to Karasmanis \cite{karasm3} pp.$148-149$ and p.$156$
section III, J. Annas \cite{annas} ch. $10,11$, and  Adams comments in \cite{thesaurus} \textit{ad.loc.}.}.
Analogous interpretation is advanced also in Plato's \textit{Tim.} $29c4$, \textit{cf}. \textit{Sophist} $232c7-9$.

There is a passage in Proclus' \textit{Commentary of the First Book of Euclid's
Elements} which seems to advance similar ideas, regarding \textsf{peras} and \textsf{apeiron} in relation to \textsf{commensurable} and
\textsf{incommensurable} magnitudes\footnote{Though that in matters of divisibility Proclus does not see the possibility
of an infinite divisibility that does not involve incommensurability, something that Plato probably
observed. For an in depth analysis of the notions of incommensurability and infinite divisibility in Proclus
see Anapolitanos and Demis \cite{anapoldem}.}, and hence it is strengthening our  conclusion ($\maltese$). Although Proclus does not refer
to Plato at all, I find it highly probable that he has in mind our passages on \textsl{peras} and
\textsl{apeiron} in the \textit{Philebus}, as we have analyzed them above.

\begin{quote}
    {\small `Mathematicals are the offspring of the limit-peratos (\textgreek{p\ea ratos}) and the unlimited-apeirian
    (\textgreek{\as peir\ia an}), but not of the primary principles alone, nor of the hidden intelligible causes, but also of secondary
    principles that proceed from them... .This is why in these orders of being there are
    ratios-logoi (\textgreek{l\oa goi}) proceeding to infinity-apeiron (\textgreek{\asa peiron}),
    but controlled by the principle of the limit-peratos. For number, beginning with unity, is capable of indefinite increase,
    yet any number you choose is finite; magnitudes-megeth\^{o}n (\textgreek{megej\wc n}) likewise are divisible without end,
    yet the magnitudes distinguished from one another are all bounded, and the actual parts of a whole
    are limited. If there were no infinity-apeirias (\textgreek{\as peir\ia as}), all magnitudes would be commensurable
    and there would be nothing inexpressible-arr\={e}ton (\textgreek{\asa rrhton}) or irrational-alogon (\textgreek{\asa logon})
    features that are thought to distinguish geometry from arithmetic.' (6. 7-22; transl. by Morrow 1970)}
\end{quote}

Proclus in his works had also systematized material from Platonic dialogues. Specifically, in his work
\textit{The Elements of Theology} ($86.16-20$, $24-26$) he clearly distinguishes the two types of infinity (\textgreek{\asa peiron})
that most probably correspond to the ones that are suggested in \textit{Philebus} and are analyzed in the sequel.
Analogous approach was advanced by Plotinus in \textit{Ennead} $2.4.15-16$.

In relation to the above, Plato in \textit{Philebus} $27e$  speaks about pleasure states that it is
\textit{apeiron} both in `\textsl{quantity and degree}', `\textgreek{ka\ig{} pl\ha jei{} ka\ig{} t\og{} m\ac llon}'.
Karasmanis in \cite{karasm} attempts to provide an explanation to this, even though he admits
that this passage is rather problematic and the possible conclusions drawn from it are not absolutely certain.
In spite of this, he  goes on providing evidence that is highly probable that \textsl{apeiron} in quantity
(\textgreek{pl\hc jos}) is something that is infinity by addition; that is, this \textsl{apeiron} is not continuous but discrete.
Furthermore, the \textsl{apeiron} in \textsl{degree} (\textgreek{t\og{} m\ac llon}) does not have definite quantities;
therefore, it is continuous or infinite by division.

Now we discuss in some extend the concepts of commensurability and incommensurability in relation to
the terms-classes of \textsl{peras} and \textsl{apeiron}.

The \textsl{incommensurability} appears as special case of of infinite divisibility. Plato in \textit{Theaet.} $147d-148e$ presents topics
from the theory of incommensurability (see Karasmanis \cite{karasm} note $16$ and Fowler \cite{fowler}).

To prove that two magnitudes are incommensurable, Greek mathematics use the technique of
anthuphairesis (\textgreek{\as njufa\ia resis}) which is used mainly to find the greatest common divisor between two numbers.
The technique of anthuphairesis is presented in Euclid's Elements as Proposition X.$2$ (see also  Karasmanis \cite{karasm} note $21$).
It is important to state that an \emph{infinite} anthuphairetic process of reciprocal subtraction between two
magnitudes shows that these magnitudes are \textsl{incommensurable}. Moreover, a \emph{finite-terminating} anthuphairetic process shows that
these magnitudes are \textsl{commensurable}\footnote{For a formal and extensive presentation of the technique of anthuphairesis
and the presentation of the proofs of the results regarding  commensurable and incommensurable magnitudes the reader
should consult Fowler \cite{fowler} ($1987$ chs $2$ and $5$), Anapolitanos \cite{anapol},
Knorr \cite{knorr} ($1975$, chs $2,4,7$), and Sinnige \cite{sinnige} ($1986$, pp. $73-80$).\label{fnt:anth}}.

It is evident that Plato apprehended these results. Specifically, the construction of incommensurable magnitudes
in \textit{Theaet.} $147d-148e$, as well as the philosophical aspect of anthuphairetic process in Plato's dialogues,
such as \textit{Sophist} $264b9-268d5$ (see Negrepontis \cite{negrep}), point to Plato's deep grasp of the technique
of anthuphairesis and his intention to apply it systematically in a general philosophical framework.
Thus, taking into account all of the above as well as the interpretation given for the phrase
`\textsl{quantity and degree}' in $27e$, one could rather safely conclude  that Plato comprehended the \emph{two types of infinity},
namely the \emph{countable and uncountable} ones\footnote{ I believe that we should view this just as Plato's attempt to discern between
two types of infinity. In modern mathematical terms, \textsl{apeiron in quantity} corresponds to  countably infinity
(with cardinality $\mathcal{N}_0$); and \textsl{apeiron in degree} corresponds to uncountably infinity
(with cardinality $2^{\mathcal{N}_0}$). For further analysis see  Karasmanis \cite{karasm}, Anapolitanos \cite{anapol}
and Anapolitanos and Demis \cite{anapoldem}.}.

In this way we see that we \emph{have two kinds of infinite divisibility} that are strongly related to magnitudes:

(1) \emph{the Zenonian infinite divisibility} that results to the Aristotelian continuum, which it points
to what we nowadays call the set of rational numbers (see Karasmanis note $4$ in \cite{karasm}).
This is a characteristic of commensurability-noting that \textsl{commensurable} magnitudes have finite anthuphairesis-and discrete infinity.
It could be considered analogous to the type of infinite that characterizes the set of rational numbers, $\mathbb{Q}$.

(2) \emph{Anthuphairetic infinite divisibility} which produces the \textsl{incommensurable} and makes the \emph{continuum dense},
thus generating  magnitudes. It could be considered analogous to the type of infinity that characterizes
the irrational numbers $\urcorner \mathbb{Q}$\footnote{where $\urcorner S$ denotes the complement of a set $S$.}.
We note that by making the continuum dense we obtain the \textsl{Real Line} $\mathbb{R}$.

Henceforth, supplementing the Aristotelian continuum-which is characterized by Zenonian infinite divisibility and corresponds
to commensurables-with the incommensurables- which are obtained via the anthuphairetic infinite divisibility-
we obtain all magnitudes. The density of the continuum and hence the construction of the \textsl{Real Line}
could be written now (using modern mathematical terminology) as $\mathbb{R}= \mathbb{Q} \cup \urcorner \mathbb{Q}$.
Recall also that $\mathbb{Q}$ has cardinality $\mathcal{N}_0$ and is a countable set, but $\urcorner\mathbb{Q}$
as well as  $\mathbb{R}$ have cardinality $2^{\mathcal{N}_0}$  and hence are uncountable sets.

Summarizing, we state the following correspondences that clarify the terms-classes
\textsl{peras}, \textsl{apeiron}:

1. \textsl{peras} $\equiv$ \textsl{commensurability} $\equiv$ \textsl{finite anthuphairetic process}.

2. \textsl{apeiron} $\equiv$ incommensurability $\equiv$ \textsl{infinite anthuphairetic process}.

Here we must note that the concepts of infinity as described above, as well as the construction of real numbers
(and in particular irrational numbers, which correspond to incommensurable magnitudes) impose a certain natural topology
on the \textsl{Real Line}\footnote{For the construction of the  \textsl{Real Line} from the set of Rational numbers
and its associate topology  we refer to Rudin \cite{rudin}. For an alternative but equivalent approach using Cauchy sequences
we recommend Moschovakis \cite{moschovakis}.}.
Within the framework of this topology the various mathematical notions, especially the concepts of density and convergence,
should be understood.
This  framework provides the appropriate meaning of the various mathematical notions involved in the formulation of the main Claim
and additionally it is the only possible one consistent with the Mathematics of Plato's Academy.

Furthermore, the concept of convergence (as apprehended in this Section) was used to argue in Section $3$B that the logics
of the TMA and the TBA arguments are not contradictory but they are rather consistent.

Of course one could naturally ask: Did Plato apprehended all these results?
I believe, as the above analysis showed, that yes, it is highly probable. We have to state that
Plato indeed had a thorough knowledge of the mathematics of his era and
especially the philosophical and foundational problems of it (see Karasmanis \cite{karasm2} and Fowler \cite{fowler}).
Our analysis leads us to suppose that \emph{Plato had a strong grasp on approaching magnitudes in terms
of incommensurability rather than the Zenonian infinite divisibility}.
Additionally, and in relation to the above, it is not unreasonable to suppose that Plato comprehended the modern mathematical notion
of \textsl{density} and through it the approach of the \textsl{continuum} as the \textsl{closure}
(in the  topology described earlier) of the \textsl{rational} numbers in \textsl{Real} numbers,
$\overline{\mathbb{Q}}=\mathbb{R}$\footnote{In Plato's thought, that associates the incommensurable magnitudes to the irrational numbers
(that have infinite anthuphairesis) and the commensurable magnitudes with rational numbers (that have finite anthuphairesis)
and in his approach in apprehending the \textsl{continuum}, we are able to discern the fundamental concepts
of the proof of the following theorem:
\begin{quote}
`Every element of the \textsl{Real Line}, and in particular every irrational number,
is the limit of a sequence of rational numbers. That is, $\forall x \in \mathbb{R}$, there is a sequence of rational numbers
$\{r_i\}_{i=0}^{\infty}$, such that $\underset{i \to +\infty} \lim r_i=x$.'
\end{quote}
For a formal proof of the theorem see Rudin \cite{rudin}.}.

\section{Conclusion}

In this section, summarizing our results, we state the main conclusions of our article.

The increasing infinite sequence of Forms $\{F_i\}_{i=0}^\infty$  constructed by  repetitive applications of the TMA
shall be understood as a \emph{convergent} one. Precisely, it converges  to the \emph{unique} $F$-Form for a particular predicate in the
framework of the topology developed in Section $4$. That is $\underset{i \to +\infty} \lim F_i=F$. In a philosophical level, we could claim
that this limit $F$-Form  should be apprehended as the constant-unchangeable  `\emph{Is}', `\textgreek{e\isc nai}'.
As we have seen in Sections $3,4$,  the convergence of the sequence $\{F_i\}_{i=0}^\infty$ to $F$ could apprehended as corresponding
to the procedure called \textsf{genesis eis ousian}- `\textgreek{g\ea nesis{} e\is s{} o\us s\ia an}', applied for each predication-Form.

Furthermore, we studied the apparent contradictory aspects of:

1. The appearances of a plurality of Forms (for each predication) that is generated by repetitive applications
of the TMA (adopting Cohen's, \cite{cohen1} thesis).

2. The intention of Plato for assuming a unique Form per Predicate in supporting the Uniqueness thesis
(as this is entailed in the TBA \textit{et.al}).

In this direction, we provided  adequate arguments defending the thesis that the logics of the TMA and the TBA are not contradictory,
but rather consistent and mutually complementary. In doing this, we have also supported that the final-limit $F$-Form
is not committed to NI axiom. Thus, we showed that the \emph{Uniqueness thesis} is safeguarded and the \emph{Abstraction Principle}
in logic is retained.

 \end{document}